\documentclass[12pt,reqno,a4wide]{amsart}

\usepackage{tikz} 

\usetikzlibrary{decorations.pathreplacing}
\usetikzlibrary{arrows,decorations.pathmorphing,snakes}
\usetikzlibrary{fit}
\usepackage{calrsfs}
\usepackage[charter]{mathdesign}

\allowdisplaybreaks

\title{Retakh's Motzkin paths and some combinatorial comments}

\author{Helmut Prodinger}
\address{Helmut Prodinger\\
	Department of Mathematical Sciences\\
	Stellenbosch University\\
	7602 Stellenbosch,	South Africa, and\\
	 Department of Mathematics and Mathematical Statistics\\
	 Umea University\\
	 907 36 Umea, 	 Sweden  }
\email{hproding@sun.ac.za}

\keywords{Motzkin numbers, Dyck paths, peaks, planar trees}
\subjclass[2020]{05A15, 05A16}

\begin{document}
	
	\begin{abstract}
		Dyck paths where peaks are only allowed on level 1 and on even-indexed levels, were introduced by Retakh
		and analysed by Zeilberger, with assistance from Ekhad. We add some combinatorial comments to the enumeration,
		which involves Motzkin numbers, in particular, about the average height of such objects.
	\end{abstract}
	
	\maketitle

	\section{Introduction}

V.~Retakh~\cite{EZ} introduced the following restricted class of Dyck paths: Peaks are only allowed on level 1 and on even-numbered levels. Here is an example, and the corresponding plane tree using the standard bijection. 
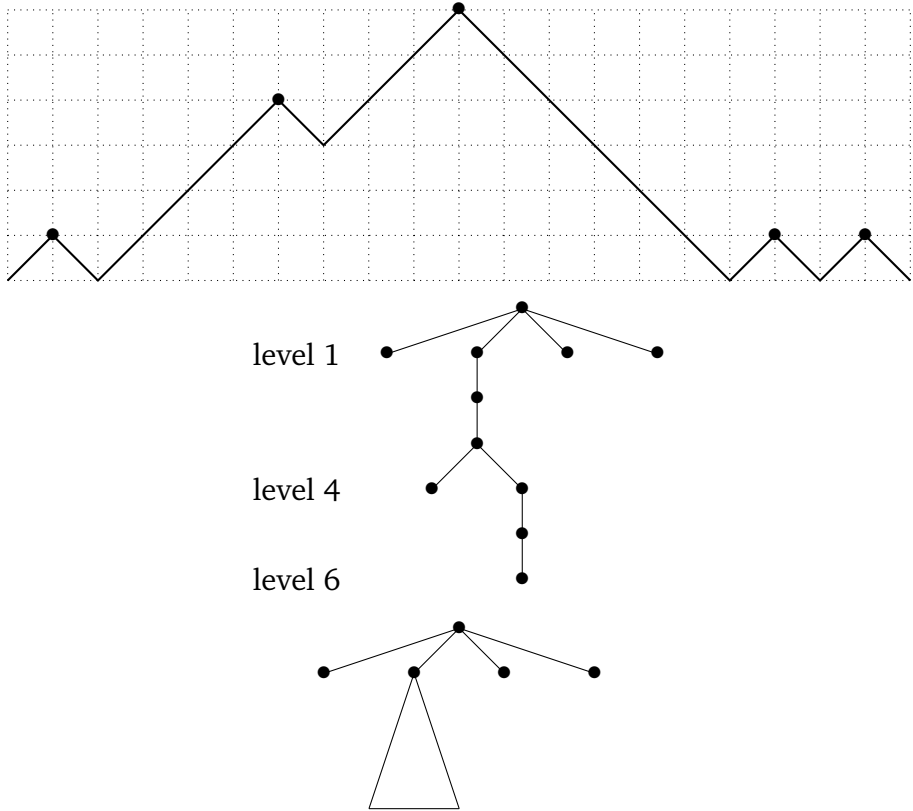
\begin{center}
	\begin{tikzpicture}[scale=0.6]
	
	\draw[step=1.cm,black,dotted] (-0.0,-0.0) grid (20.0,6.0);


	\draw[thick] (0,0) -- (1,1) -- (2,0)-- (3,1)-- (4,2)-- (5,3)-- (6,4)-- (7,3)-- (8,4)-- (9,5)-- (10,6)-- (11,5)-- (12,4)-- (13,3)-- (14,2)-- (15,1)-- (16,0)-- (17,1)-- (18,0)-- (19,1)-- (20,0);
	
	\node at (1,1){$\bullet$};
	\node at (6,4){$\bullet$};
	\node at (10,6){$\bullet$};
	\node at (17,1){$\bullet$};
		\node at (19,1){$\bullet$};
	\end{tikzpicture}
	\end	{center}
	
	\begin{center}
		\begin{tikzpicture}[scale=0.6]
		\node (1) at (0, 0){$\bullet$};
		\node (2) at (-3, -1){$\bullet$};
		\node (3) at (-1, -1){$\bullet$};
		\node (4) at (1, -1){$\bullet$};
		\node (5) at (3, -1){$\bullet$};
		\node (6) at (-1, -2){$\bullet$};
		\node (7) at (-1, -3){$\bullet$};
		\node (8) at (-2, -4){$\bullet$};
			\node (9) at (0, -4){$\bullet$};	
			\node (10) at (0, -5){$\bullet$};
			\node (11) at (0, -6){$\bullet$};
			
			\node  at (-5, -1){level $1$};
\node  at (-5, -4){level $4$};	
\node  at (-5, -6){level $6$};		
			
		\draw[-] (1.center) to (2.center);
		\draw[-] (1.center) to (3.center);
		\draw[-] (1.center) to (4.center);
		\draw[-] (1.center) to (5.center);
		\draw[-] (3.center) to (6.center);
		\draw[-] (6.center) to (7.center);
		\draw[-] (8.center) to (7.center);
		\draw[-] (9.center) to (7.center);
			\draw[-] (9.center) to (10.center);
				\draw[-] (10.center) to (11.center);
		\end{tikzpicture}
		
	\end{center}
	
	\begin{center}
		\begin{tikzpicture}[scale=0.6]
		\node (1) at (0, 0){$\bullet$};
		\node (2) at (-3, -1){$\bullet$};
		\node (3) at (-1, -1){$\bullet$};
		\node (4) at (1, -1){$\bullet$};
		\node (5) at (3, -1){$\bullet$};
	
		\node (8) at (-2, -4){  };
		\node (9) at (0, -4){ };

		\draw[-] (1.center) to (2.center);
		\draw[-] (1.center) to (3.center);
		\draw[-] (1.center) to (4.center);
		\draw[-] (1.center) to (5.center);
		\draw[-] (3.center) to (8.center);
		\draw[-] (3.center) to (9.center);
		\draw[-] (8.center) to (9.center);
		\end{tikzpicture}
		
	\end{center}
	Ekhad and Zeilberger \cite{EZ} proved  recently that these restricted paths are enumerated by Motzkin numbers. 
	Recall that the generating function of the Motzkin numbers $M(z)$ according to length satisfies $M=1+zM+z^2M^2$ and thus
	\begin{equation*}
	M(z)=\frac{1-z-\sqrt{1-2z-3z^2}}{2z^2}.
	\end{equation*}
	
	In this note, I want to present a few additional observations, also including the height of the paths (or the associated plane trees).
	First, we are going to confirm the connection to Motzkin paths:	
	
Since the level 1 is somewhat special, we only consider  trees as symbolized by the triangle. We will use two generating functions, to deal with the odd/even situation. We have
\begin{equation*}
F=\frac{zG}{1-G}\quad\text{and}\quad G=\frac{z}{1-F}.
\end{equation*}
Solving, we find $F(z)=z^2M(z)$ and the total generating function is
	\begin{equation*}
	\frac z{1-z}\sum_{r\ge0}\Big(\frac{F}{1-z}\Big)^r=zM(z),
	\end{equation*}
	as predicted. Recall that the number of nodes in trees is always one more than the half-length of the corresponding Dyck path.

We will       compute the average height of such restricted paths, using singularity analysis of generating functions, as in \cite{FlOd,FS}. Whether we define the height in terms of the maximal chain of edges resp.\ nodes only makes a difference of one,
and we will only compute the average height according to the leading term of order $\sqrt n$. For readers who wish to see how more terms could be computed, at least in principle, we cite~\cite{Prodinger-ars}.
	
\section{The height}	

Now we will use the substitution $z=\frac{v}{1+v+v^2}$, which occured for the first time in \cite{Prodinger-three}, but has been used more recently in different models where Motzkin numbers are involved \cite{Prodinger-Deutsch, HHP, Baril}.  For example, we have simply $M(z)=1+v+v^2$.
We define 
\begin{equation*}
G_{k+1}=\dfrac{z}{1-\dfrac{zG_k}{1-G_k}},\quad\text{with} \quad G_1=z.
\end{equation*}
There is a simple formula, viz.
\begin{equation*}
G_k=\frac v{1+v}\frac{1-v^{2k}}{1-v^{2k+1}}.
\end{equation*}
This is easy to prove by induction, which we will do for the convenience of the reader.
The start is
\begin{equation*}
	G_1=\frac v{1+v}\frac{1-v^{2}}{1-v^{3}}=\frac v{1+v}\frac{1+v}{1+v+v^{2}}=\frac{v}{1+v+v^{2}}=z.
\end{equation*}
And now
\begin{align*}
G_{k+1}&=\dfrac{z}{1-\dfrac{zG_k}{1-G_k}}=\dfrac{z(1-G_k)}{1-(1+z)G_k}=\frac{v}{1+v+v^{2}}\dfrac{1-\frac v{1+v}\frac{1-v^{2k}}{1-v^{2k+1}}}{1-\frac{(1+v)^2}{1+v+v^2}\frac v{1+v}\frac{1-v^{2k}}{1-v^{2k+1}}}\\
&=\frac{v}{1+v}\dfrac{1+v-v\frac{1-v^{2k}}{1-v^{2k+1}}}{1+v+v^{2}-v(1+v)\frac{1-v^{2k}}{1-v^{2k+1}}}
=\frac{v}{1+v}\dfrac{(1+v)(1-v^{2k+1})-v(1-v^{2k})}{(1+v+v^{2})(1-v^{2k+1})-v(1+v)(1-v^{2k})}\\
&=\frac{v}{1+v}\dfrac{1-v^{2k+2}}{1-v^{2k+3}},
\end{align*}
as claimed. From this we also get
\begin{equation*}
F_k=\frac{zG_k}{1-G_k}=\frac{v^2}{1+v+v^2}\frac{1-v^{2k}	}{1-v^{2k+2}}.
\end{equation*}
For $k\ge1$, $F_k$ is the generating function of trees (like in the triangle) of height $\le 2k$. 

Note that the height is currently counted in terms of nodes;
\begin{equation*}
F_1=\frac{z^2}{1-z},
\end{equation*} 
which describes a root with $\ell\ge1$ leaves attached to the root.

Now we incorporate the irregular beginning of the tree and compute
\begin{equation*}
\frac z{1-z}\sum_{r\ge0}\Big(\frac{F_h}{1-z}\Big)^r=\frac z{1-z}\dfrac1{1-\dfrac{F_h}{1-z}}=v\frac{1-v^{2h+2}}{1-v^{2h+4}}.
\end{equation*}
From here onwards it seems to be more natural to define the height of the whole tree in terms of the number of \emph{edges}, 
and then the quantity we just derived is the generating function of all trees with height $\le 2h$, for $h\ge1$. Note that the limit
$h\to\infty$ gives us simply $v=zM(z)$, which is consistent. There is also a contribution of trees of height $\le 1$, namely 
$\frac{z}{1-z}=\frac{v}{1+v^2}$, but this term is, when we compute the average height, irrelevant and only contributes to the error term, as we only compute  the leading term, which is of order $\sqrt n$.

So, apart from normalization, we are led to investigate
\begin{align*}
	\sum_{h\ge1}&2h\bigg[v\frac{1-v^{2h+2}}{1-v^{2h+4}}-v\frac{1-v^{2h}}{1-v^{2h+2}}\bigg]
=2v(1-v^{-2})\sum_{h\ge1}h\bigg[\frac{v^{2h+4}}{1-v^{2h+4}}-\frac{v^{2h+2}}{1-v^{2h+2}}\bigg]\\
&=2v(1-v^{-2})\sum_{h\ge0}h\frac{v^{2h+4}}{1-v^{2h+4}}
-2v(1-v^{-2})\sum_{h\ge0}(h+1)\frac{v^{2h+4}}{1-v^{2h+4}}\\
&=-2v+\frac{2(1-v^{2})}{v}\sum_{h\ge1}\frac{v^{2h}}{1-v^{2h}}.
\end{align*}

Note that we could get explicit coefficients fr0m here, using trinomial coefficient, $\binom{n,3}{k}=[v^k](1+v+v^2)^n$ (notation from \cite{Comtet-book}). To show the reader how this works, we compute
\begin{align*}
[z^{n+1}]&\frac{1-v^{2}}{v}\sum_{h\ge1}\frac{v^{2h}}{1-v^{2h}}
=\frac1{2\pi i}\oint\frac{dz}{z^{n+2}}\frac{1-v^{2}}{v}\sum_{h\ge1}\frac{v^{2h}}{1-v^{2h}}\\
&=\frac1{2\pi i}\oint dv(1-v^2)^2\frac{(1+v+v^2)^{n}}{v^{n+3}}\sum_{h\ge1}\sum_{h,k\ge1}v^{2hk}\\
&=[v^{n+2}](1-2v^2+v^4)\sum_{h\ge1}d(h)v^{2h}(1+v+v^2)^{n}\\
&=\sum_{h\ge1}d(h)\bigg[\binom{n,3}{n+2-2h}-2\binom{n,3}{n-2h}+\binom{n,3}{n-2-2h}\bigg].
\end{align*}
Note that $d(h)$ is the number of divisors of $h$. We will, however, not use this explicit form.
The expression as derived before,
\begin{equation*}
	-2v+\frac{2(1-v^{2})}{v}\sum_{h\ge1}\frac{v^{2h}}{1-v^{2h}},
\end{equation*}
has to be expanded around $v=1$, which is a standard application of the Mellin transform. Details are worked out in
\cite{HPW}, for example:
\begin{equation*}
\sum_{h\ge1}\frac{v^{2h}}{1-v^{2h}}=\sum_{k\ge1}d(k)v^{2k}\sim-\frac{\log(1-v^2)}{1-v^2}
\sim-\frac{\log(1-v)}{2(1-v)}.
\end{equation*}
Note again that $d(k)$ is the number of divisors of $k$. Consequently
\begin{equation*}
-2v+\frac{2(1-v^{2})}{v}\sum_{h\ge1}\frac{v^{2h}}{1-v^{2h}}\sim
-2\log(1-v).
\end{equation*}
We have $1-v\sim \sqrt 3\sqrt{1-3z}$, and $z=\frac13$ is the relevant singularity when discussing Motzkin numbers. We can continue
\begin{equation*}
-2v+\frac{2(1-v^{2})}{v}\sum_{h\ge1}\frac{v^{2h}}{1-v^{2h}}\sim -\log(1-3z).
\end{equation*}
The coefficient of $z^n$ in this is $\frac{3^n}{n}$. This has to be divided by
\begin{equation*}
[z^n]zM(z)=[z^{n-1}]M(z)\sim \frac{3^{n+\frac12}}{2\sqrt{\pi }n^{3/2}},
\end{equation*}
with the final result for the average height of the restricted Dyck paths (\`a la Retakh):
\begin{equation*}
\sim 2\sqrt{\frac{\pi n}{3}}.
\end{equation*}
Recall \cite{Prodinger-three} that the average height of Motzkin paths of length $n$ is asymptotic to
\begin{equation*}
\sqrt{\frac{\pi n}{3}}.
\end{equation*}

\section{The number of leaves}
We can use a second variable, $u$, to count the number of leaves. Then we have
\begin{equation*}
F(z,u)=\frac{zG(z,u)}{1-G(z,u)}\quad\text{and}\quad G(z,u)=zu+\frac{zF(z,u)}{1-F(z,u)},
\end{equation*}
which leads to
\begin{equation*}
	F(z,u)= {\frac {1-zu-{z}^{2}+{z}^{2}u-\sqrt {1-2zu-2{z}^{2}-2{z}^{
					2}u+{z}^{2}{u}^{2}-2{z}^{3}u+2{z}^{3}{u}^{2}+{z}^{4}-2{z}^{4}u+{
					z}^{4}{u}^{2}}}{2(1-zu+z)}}.
		\end{equation*}
Bringing the irregular beginning also into the game leads to
\begin{align*}
\frac{z}{1-zu}\sum_{r\ge0}\Big(\frac{F}{1-zu}\Big)^r+zu-z.
\end{align*}
This is an ugly expression that we do not display here. But we can compute the average number of leaves, by differentiation w.r.t. $u$, followed by setting $u=1$:
\begin{equation*}
R:={\frac {v \left( 1+v \right)   ( 1-v+2{v}^{2}-{v}^{3}) }{
		(1-v)  \left( 1+v+{v}^{2} \right) }}.
\end{equation*}
The coefficient of $z^n$ in this can be expressed in terms of trinomial coefficients, if needed. But we only compute an asymptotic formula, to keep this section short. Expanding around $v=1$, we find
\begin{equation*}
R\sim\frac23\frac{1}{1-v}\sim\frac23\frac{1}{\sqrt3\sqrt{1-3z}},
\end{equation*}
and thus
\begin{equation*}
[z^n]R\sim\frac23\frac{1}{\sqrt3}3^n\frac1{\sqrt{\pi n}}.
\end{equation*}
We divide this again by
\begin{equation*}
\frac{3^{n+\frac12}}{2\sqrt{\pi} n^{3/2}}
\end{equation*}
with the result
\begin{equation*}
\frac49n,
\end{equation*}
which is the asymptotic number of leaves in a Retakh tree of size $n$. Recall that for unrestricted planar trees, the result is
$\frac n2$, which is a folklore result using Narayana numbers. So the constant in the restricted case, $\frac49$, is a bit smaller than 
$\frac12$.

With some effort, more precise approximations could be obtained, as well as the variance. This might be a good project for a student.

\section{Conclusion}

I am always happy to see a new occurrence of Motzkin numbers and that the methods that I learnt more than 40 years ago from Knuth, Flajolet, and others still work.


\end{document}